\documentclass[letterpaper,12pt]{amsart}
\usepackage{epic,eepic,amsfonts,color,latexsym,mathrsfs,amssymb,amsmath,epsfig,rotating,hyperref}
\newtheorem{thm}{Theorem}[section]
\newtheorem{dfn}[thm]{Definition}

\newtheorem{cor}[thm]{Corollary} 

\newtheorem{prop}[thm]{Proposition}
\newtheorem{lemma}[thm]{Lemma}
\newtheorem{ex}[thm]{Example}
\newtheorem{forthm}{Forsg\aa{}rd's Theorem} 
  
\newtheorem{tame}{Semi-Algebraic Tameness Theorem} 
  
\newtheorem{sop}{Proposition BCS}

\makeatletter
\renewcommand*\env@matrix[1][c]{\hskip -\arraycolsep
  \let\@ifnextchar\new@ifnextchar
  \array{*\c@MaxMatrixCols #1}}
\makeatother

\definecolor{red}{rgb}{.5,0,0} 
\definecolor{green}{rgb}{0,.4,0} 
\definecolor{blue}{rgb}{0,0,.5} 
\newcommand{\um}{u_{\mathrm{min}}}
\newcommand{\up}{u_{\mathrm{max}}}
\newcommand{\gm}{\gamma_{\mathrm{inf}}}
\newcommand{\gp}{\gamma_{\mathrm{sup}}}

\newcommand{\thth}{^{\text{\underline{th}}}}

\newcommand{\np}{{\mathbf{NP}}}

\newcommand{\eps}{\varepsilon}

\newcommand{\Log}{\mathrm{Log}}

\newcommand{\Q}{\mathbb{Q}}
\newcommand{\R}{\mathbb{R}}
\newcommand{\g}{\tilde{g}}

\newcommand{\C}{\mathbb{C}}
\newcommand{\N}{\mathbb{N}}
\newcommand{\Z}{\mathbb{Z}}
\newcommand{\bO}{\mathbf{O}}
\newcommand{\cA}{\mathcal{A}}

\newcommand{\Zn}{\Z^n}

\newcommand{\Rn}{\R^n}

\newcommand{\Cn}{\C^n}

\newcommand{\Cs}{\C^*}

\newcommand{\cN}{{\mathcal{N}}}

\newcommand{\re}{{\mathbf{Re}}}
\newcommand{\im}{{\mathbf{Im}}}

\renewcommand{\qed}{$\blacksquare$}
\newcommand{\dia}{$\diamond$}
\newcommand{\amoeba}{\mathrm{Amoeba}}

\newcommand{\trop}{\mathrm{Trop}}


\author{Alperen A. Erg\"ur}
\address{Technische Universit\"at Berlin, Institut f\"ur Mathematik,
Sekretariat MA 3-2, Strasse des 17. Juni 136 10623 Berlin, Germany} 
\email{erguer@math.tu-berlin.de}
\author{Grigoris Paouris}
\address{Department of Mathematics,
Texas A\&M University TAMU 3368,
College Station, Texas \ 77843-3368, USA.}  
\email{grigoris@math.tamu.edu}
\author{J.\ Maurice Rojas}
\address{Department of Mathematics,
Texas A\&M University TAMU 3368,
College Station, Texas \ 77843-3368, USA.}  
\email{rojas@math.tamu.edu}
\thanks{A.E.\ was partially supported by NSF grant CCF-1409020, and Einstein Foundation, Berlin. J.M.R.\ was partially supported by NSF grant CCF-1409020, and LABEX MILYON (ANR-10-LABX-0070) of 
Universit\'e de Lyon, within the program ``Investissements d'Avenir'' (ANR-11-IDEX-0007) operated by ANR. G.P.\ was partially supported by BSF grant 2010288 and NSF CAREER grant DMS-1151711.}

\title[Tropical Varieties for Exponential Sums]{\mbox{}\\
\vspace{-1.3in}
Tropical Varieties for Exponential Sums\\ 
} 

\dedicatory{In memory of Joel Zinn (March 16, 1946 -- December 5, 2018), 
beloved friend and brilliant colleague.} 

\date{\today{} version of paper published in Math. Ann. in 2020: DOI 
10.1007/s00208-019-01808-5 . } 

\begin{document}

\begin{abstract} 
We study the complexity of approximating complex zero sets of certain 
$n$-variate exponential sums. We show that the real part,  
$R$, of such a zero set can be\linebreak 
approximated by the $(n-1)$-dimensional 
skeleton, $T$, of a polyhedral subdivision of $\mathbb{R}^n$. In particular, 
we give an explicit upper bound on the Hausdorff distance: $\Delta(R,T)
\!=\!O\!\left(t^{3.5}/\delta\right)$, where $t$ and $\delta$ are respectively 
the number of terms and the minimal spacing of the frequencies of $g$. 
On the side of computational complexity, we show that even the 
$n\!=\!2$ case of the membership problem for $R$ is undecidable
in the Blum-Shub-Smale model over $\R$, whereas membership and distance 
queries for our polyhedral approximation $T$ can be decided in polynomial-time 
for any fixed $n$.  
\end{abstract} 

\maketitle

\vspace{-1.3cm} 
\section{Introduction} 
We study zero sets of exponential sums of the form 
$g(z):=\sum_{j=1}^t e^{a_j\cdot z + \beta_j}$ 
where $z\!\in\!\mathbb{C}^n$, $a_j \in \mathbb{R}^n$, the $a_j$ are pair-wise 
distinct, $\beta_j \in \mathbb{C}$, 
and $a_j\cdot z$ denotes the usual Euclidean inner product in 
$\mathbb{C}^n$. We call $g$ an {\em $n$-variate exponential $t$-sum}, 
$a_j$ a {\em frequency} of $g$, $\{a_1,\ldots,a_t\}$ the {\em spectrum} of $g$, 
and $\delta(g)\!:=\!\min_{p\neq q} |a_p-a_q|$ the 
{\em minimal spacing of the frequencies of $g$}. (Throughout this paper, 
we use $|\cdot|$ for the standard $\ell_2$-norm on $\C^N$ for any 
$N\!\in\!\N$.)  
We also call the $\beta_j$ the {\em coefficients} of $g$. 
One can think of $g$ as an analogue of a polynomial with real exponents, and 
hope to use algebraic intuition to derive new metric results in the broader 
setting of exponential sums. We shall do so by combining results on 
random projections with some new extensions of classical univariate polynomial 
bounds. 

Exponential sums appear across pure and applied mathematics. For instance, 
exponential sums (in the form above) occur in the 
calculation of $3$-manifold invariants (see, e.g., \cite[Appendix A]{mcmullen} 
and \cite{hadari}), and have been studied from the point of view of 
Diophantine Geometry, Model Theory, and Computational Algebra, 
(see, e.g., \cite{richardson,macwil,wilkie,zilberolder,achatz,zilber,scanlon,
habegger}). Also, the non-lattice Dirichlet polynomials appearing in the 
study of fractal strings \cite{lv} are a special case of the exponential 
sums we consider here. An application to radar antennae \cite{forsythe,hwang} 
---  finding the directions of a set of unknown signals --- reduces to  
finding the zeroes of a univariate exponential sum, with 
frequencies depending on the location of the sensors of the antenna. 
Approximating roots of multivariate exponential sums is 
also a fundamental computational problem in Geometric Programming 
\cite{duffin,chiang,boyd}.  

For any analytic function $g$ on $\Cn$ let $Z(g)$ denote the 
set of complex zeroes of $g$. Also,\linebreak 
\scalebox{.95}[1]{for any $W\!\subseteq\!\Cn$, we define its {\em real part} 
to be  $\re(W)\!:=\!\{(\re(z_1),\ldots,\re(z_n))\; | \; 
(z_1,\ldots,z_n)\!\in\!W\}$.}\linebreak  
One can wonder if exact computation with the roots of 
exponential sums is possible using only field operations and comparisons 
over $\R$, or if approximation is truly necessary. 
Exact computation turns out to be intractable, relative to a standard  
computational model (the {\em BSS model over $\R$} \cite{bcss}), 
already in the special case of two variables and three terms. 
\begin{thm} \label{thm:mem} 
Determining, for arbitrary input $r_1,r_2\!\in\!\R$, whether 
$(r_1,r_2)$ lies in\linebreak 
$\re(Z(1-e^{z_1}-e^{z_2}))$ is undecidable\footnote{\cite{poonen} provides an
excellent survey on undecidability, in the classical Turing model, geared
toward non-experts in complexity theory.} 
in the BSS model over $\R$. 
\end{thm}

\noindent 
We prove Theorem \ref{thm:mem} in Section \ref{sub:mem}. There 
are certainly tractable special cases of the preceding problem, such 
as when the $r_i\!=\!\log s_i$ for some positive rational $s_i$ 
(see, e.g., \cite{theobald,wolfsos} and \cite[Thm.\ 1.9]{aknr}). 
Similarly, the famous Lindemann-Weierstrass Theorem tells us that 
$e^{r_1}+e^{r_2}$ is transcendental when $r_1,r_2\!\in\!\R$  
are distinct and algebraic. However, checking whether 
$e^{r_1}+e^{r_2}\!\stackrel{?}{\in}\!\Q(r_1,r_2)$ for {\em arbitrary}  
distinct transcendental $r_1,r_2\!\in\!\R$ --- using only {\em finitely} many 
rational operations and inequality checks in $\Q(r_1,r_2)$ --- is already an 
open question. Theorem \ref{thm:mem} thus highlights the need for approximation 
if one wants to work with roots of exponential sums in complete generality.  

A natural question then is whether one can {\em efficiently} approximate the 
zero set of an exponential sum. For instance, can we at least decide --- 
perhaps within polynomial-time --- whether a given point is 
close to the real part of the zero set of an exponential sum? Our main 
algorithmic and quantitative results (Theorems \ref{thm:cxity2} and 
\ref{thm:intro2}) show that this is indeed the case, at least in a coarse 
sense: We derive a polyhedral structure that can be considered as a 
first-order approximation to the real part of the zero set, so that 
higher-order numerical iterative methods can be deployed when higher precision 
is needed in a specific application. 

Clearly, $Z(g)$ is empty when $t\!=\!1$.  
That polyhedra arise from the real parts of zero sets of exponential 
sums is most easily seen in the special case of $t\!=\!2$ terms: Since 
$\left|\pm e^\beta\right|\!=\!e^{\re(\beta)}$, the equality 
$e^{a_1\cdot z+\beta_1}+e^{a_2\cdot z+\beta_2}\!=\!0$ implies 
$e^{a_1\cdot \re(z)+\re(\beta_1)}\!=\!e^{a_2\cdot \re(z)+\re(\beta_2)}$, and 
we thus obtain the following basic fact after taking logarithms: 
\begin{prop}  
\label{prop:aff} 
If $g(z)\!=\!e^{a_1\cdot z+\beta_1}+e^{a_2\cdot z+\beta_2}$ for some distinct 
$a_1,a_2\!\in\!\Rn$, and $\beta_1,\beta_2\!\in\!\C$, then 
$\re(Z(g))$ is the affine hyperplane  
$\left\{\left. u\!\in\!\Rn\; \right| \; 
(a_1-a_2)\cdot u+\re(\beta_1-\beta_2)\!=\!0\right\}$. \qed 
\end{prop} 

Before stating our main metric results in arbitrary dimension, it will be 
useful to observe some of the intricacies present already in the univariate 
case. 

\subsection{Clustering of Real Parts in One Variable} 
The simple sum $e^{z_1}-1$ shows that the imaginary part $\im(Z(g))$ can be 
infinite already in the univariate case, unlike the polynomial setting. A more 
subtle phenomenon, however, is that $\re(Z(g))$ need not even be closed.  
\begin{prop} 
\label{prop:dumb} 
$X\!:=\!\re\!\left(Z\!\left(e^{\sqrt{2}z_1}+e^{\sqrt{3} z_1}+
e^{\sqrt{5}z_1}\right)\right)$ 
is countably infinite, contained in the open interval 
$\left(-\frac{\log 2}{\sqrt{3}-\sqrt{2}},\frac{\log 2}{\sqrt{3}-\sqrt{2}}
\right)$ ($\subset\!(-2.181,2.181)$), and dense in the open interval 
$(-1.06,1.06)$. In particular, $X$ does not contain all its limit points.  
\end{prop} 

\noindent 
We prove Proposition \ref{prop:dumb} in Section \ref{sec:uni}. Another 
subtlety behind $\re(Z(g))$ is that finding its points in the special  
case where $n\!=\!1$ and the spectrum of $g$ lies in $\Z$ 
is the same as finding the 
logarithms of the absolute values of the complex roots of a polynomial. In 
particular, just deciding  $0\stackrel{?}{\in}\re(Z(g))$ in this special case 
is already $\np$-hard \cite{plaisted}.   

A natural trick we will soon justify is that we can predict real parts 
by examining pairs of terms of $g$ with large absolute value, in order 
to locally reduce to the two-term case:    
\begin{dfn} 
Let us define, for any $n$-variate exponential $t$-sum $g$, with $t\!\geq\!2$, 
its\linebreak 
{\em tropical variety} as\\ 
\mbox{}\hfill $\trop(g)\!:=\!\re\left(\left\{z\!\in\!\Cn\; : \;   
\max_j\left|e^{a_j\cdot z+\beta_j}\right| \text{ is attained at at least
two distinct } j\right\}\right)$.\hfill\dia 
\end{dfn} 

\noindent 
The calculation preceding Proposition \ref{prop:aff} in fact yields 
$\trop(g)\!=\!\re(Z(g))$ when $t\!=\!2$. 

\noindent 
\begin{minipage}[b]{0.65\linewidth}
\vspace{0pt}
\mbox{}\hspace{.5cm}More generally, among many other equivalent 
characterizations, 
$\trop(g)$ can also be defined as the set of points at which the 
piece-wise linear function $\cN_g : \Rn\longrightarrow \R$ defined by 
$\cN_g(u)\!:=\!\max_j\left\{a_j\cdot u + \re(\beta_j)\right\}$ is 
non-differentiable. So, for $n\!=\!1$, the graph of $\cN_g$ is concave upward, 
with at most $t-1$ ``corners,'' and thus $\trop(g)$ consists of at most $t-1$ 
points. For instance, $g(z_1)\!:=\!(e^{z_1}+1)^2$ implies that 
$\cN_g(u)\!=\!\max\{0,u+\log 2,2u\}$ (with graph drawn to the right) and 
thus $\trop(g)\!=\!\{\pm \log 2\}$.
\end{minipage}\hspace{.3cm}   
\begin{minipage}[b]{0.3\linewidth}
\vspace{0pt}
\epsfig{file=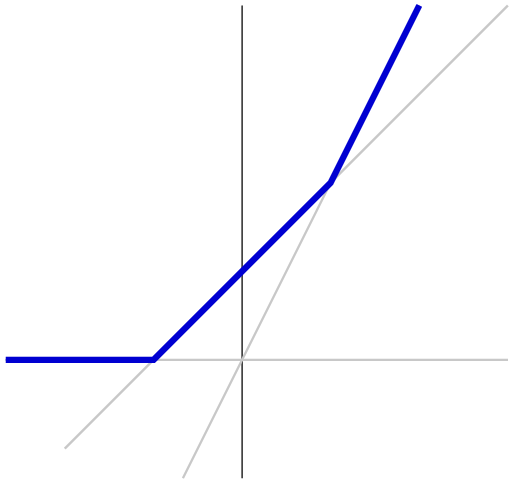,height=1.7in,clip=} 
\end{minipage} 

Computing $\trop(g)$ when $n\!=\!1$ is thus no harder than computing 
a convex hull in $\R^2$, and 
$\re(Z(g))$ turns out to always accumulate predictably near $\trop(g)$. 
In what follows, we use $\#S$ for the cardinality of a set $S$. 
\begin{thm} 
\label{thm:uni} 
Suppose $g$ is any univariate exponential $t$-sum with spectrum 
$\{a_1,\ldots,a_t\}\!\subset\!\R$, minimal frequency spacing 
$\delta(g)\!:=\!\min_{p\neq q} |a_p-a_q|$, and  
$t\!\geq\!3$. Let $s\!:=\!\#\trop(g)$,\linebreak  
$\um\!:=\!\min\trop(g)$, $\up\!:=\!\max\trop(g)$, and let $U_g$ be the union 
of open intervals\\  
\mbox{}\hfill $\left(\um-\frac{\log 2}{\delta(g)},\up+\frac{\log 2}{\delta(g)}
\right) \cap \bigcup\limits_{u\in \trop(g)} \!\!\!\!\!\!\!\left(u-\frac{\log 3}
{\delta(g)}, u+\frac{\log 3}{\delta(g)}\right)$.\hfill\mbox{}\\  
Then $1\!\leq\!s\!\leq\!t-1$ and:\\   
\mbox{}\hspace{1.5cm}(1) $\re(Z(g))\!\subset\!U_g$.\\ 
\mbox{}\hspace{1.5cm}(2) $\re(Z(g))$ has at least one point in each connected 
component of $U_g$.\\  
\mbox{}\hspace{1.5cm}(3) For any $u\!\in\!\trop(g)$ there is a root 
$\zeta\!\in\!\C$ of $g$ with\\ 
\mbox{}\hfill $|u-\re(\zeta)|<\frac{(\log 9)s-\log\frac{9}{2}}{\delta(g)}\leq 
\frac{(\log 9)t-\log\frac{81}{2}}{\delta(g)}<(2.2t-3.7)/\delta(g)$.\hfill\mbox{}  
\end{thm} 

\noindent 
We prove Theorem \ref{thm:uni} in Section \ref{sub:uni}. 
The constants $\log 2$ and $\log 3$ in the definition of the neighborhood 
$U_g$ above are in fact optimal:  
\begin{lemma} 
\label{lemma:opt} (See, e.g., \cite[Cor.\ 2.3(c) \& 
Lemma 2.5]{aknr}.) Consider any real $\delta\!>\!0$, any integer 
$t\!\geq\!2$, and the exponential sums\\  
\mbox{}\hspace{.85cm}$g_{1,t}(z_1) \; \; := \; \; e^{(t-1)\delta z_1} 
 - e^{(t-2)\delta z_1} 
-\cdots - e^0$ \ \ , \ \ $g_{2,t}(z_1) := g_{1,t}(-z_1)$ \ \ , \ \ and   
\begin{eqnarray*} 
g_{3,t}(z_1) & := & 1+e^{\delta z_1}+\cdots+e^{(t-1)\delta z_1}
-e^{t\delta z_1}
+e^{(t+1)\delta z_1-1\cdot \log 9}+\cdots+e^{(t+t)\delta z_1-t\log 9}. 
\end{eqnarray*} 
Then we have:\\ 
(1) $\trop(g_{1,t})\!=\!\trop(g_{2,t})\!=\!\{0\}$ but 
$\re(Z(g_{1,t}))$ (resp. $\re(Z(g_{2,t}))$) contains points strictly\\
\mbox{}\hspace{.7cm}increasing (resp.\ strictly 
decreasing) toward a limit of $\frac{\log 2}{\delta}$ (resp.\ 
$-\frac{\log 2}{\delta}$) as $t\longrightarrow \infty$.\\  
(2) $\trop(g_{3,t})\!=\!\left\{0,\frac{\log 9}{\delta}
\right\}$ and $\re(Z(g_{3,t}))\cap\left[-\frac{\log 3}{\delta},\frac{\log 3}
{\delta}\right]$ is empty. However, for any $\eps\!>\!0$,\\ 
\mbox{}\hspace{.7cm}there is a $t\!\in\!\N$ 
such that $\re(Z(g_{3,t}))\cap\left(\frac{\log 3}{\delta}-\eps,
\frac{\log 3}{\delta}+\eps\right)$ is non-empty. \qed 
\end{lemma}   

\noindent 
We note that \cite[Cor.\ 2.3(c) \& Lemma 2.5]{aknr}, while 
phrased in terms of univariate polynomials $f(x_1)$, directly yield  
Assertions (1) and (2) above upon substituting $x_1\!=\!e^{\delta z_1}$. 

The clustering of $\re(Z(g))$ about $\trop(g)$ persists in higher dimension.  

\subsection{Efficiently Finding Clusters of Real Parts in Arbitrary Dimension} 
Our definition of tropical variety generalizes an 
earlier version defined just for polynomials: When the spectrum 
of $g$ lies in $\Zn$, one can associate to 
our exponential sum $g$ the Laurent 
polynomial $f(x)\!:=\!\sum^t_{j=1} e^{\beta_j} 
x^{a_j}\!\in\!\C\!\left[x^{\pm 1}_1,\ldots,x^{\pm 1}_n\right]$. 
Recall that the {\em amoeba} of $f$ is the set\\ 
\mbox{}\hfill  
$\amoeba(f)\!:=\!\{(\log|x_1|,\ldots,\log|x_n|)\; | \; f(x_1,\ldots,x_n)\!=\!0; 
x_1,\ldots,x_n\!\in\!\Cs\}$. \hfill\mbox{}\\  
It is then clear that, under these restrictions,  
$\re(Z(g))\!=\!\amoeba(f)$. {\em Tropical geometry} 
(see, e.g., \cite{virologpaper,passarerullgard,ekl,payne,tropical1,bakerrumely,
barba,macsturmf,aknr}) enables algebraic varieties over various complete 
algebraically closed fields (such as $\C$, $\C\langle\langle t\rangle\rangle$, 
or $\C_p$, to name a few) to be approached polyhedrally. In our notation 
here, defining $\trop(f)\!:=\!\trop(g)$ results in the {\em 
Archimedean} tropical variety of $f$, whose metric aspects were studied 
in \cite{aknr}. This kind of tropical variety over $\C$ can be 
traced back to 1893 work of Hadamard revealing how to 
polyhedrally approximate products of 
norms of complex roots of univariate polynomials \cite{hadamard}. 

Our $\trop(\cdot)$ here is thus a small first step toward extending   
tropical methods from polynomial functions to certain exponential sums.  
It should be noted that the theory of $\cA$-discriminants \cite{gkz94} 
now also has a generalization to exponential sums \cite{rojasrusek}, 
and these generalizations have led to sharper bounds in 
real fewnomial theory \cite{fnr}.  

Recall that the {\em affine span} of a point set $S\!\subset\!\Rn$ 
is the smallest affine subspace of $\Rn$ containing $S$. 
Via polyhedral duality (see, e.g., 
\cite{grunbaum,ziegler,triang}), an immediate consequence of our 
characterization of $\trop(g)$ via the graph of $\cN_g$ is the following fact:  
\begin{prop} 
Let $d$ be the dimension of the affine span of the spectrum of a 
real $n$-variate exponential $t$-sum $g$. Then  
$\trop(g)$ is a polyhedral complex of 
pure dimension $n-1$, and is connected when $d\!\geq\!2$. \qed 
\end{prop} 

\begin{dfn} For any $n$-variate exponential $t$-sum $g$, let
$\Sigma(\trop(g))$ denote the polyhedral complex whose cells are
exactly the (possibly improper) faces of the closures of the
connected components of $\Rn\!\setminus\!\trop(g)$. \dia
\end{dfn} 

We can now make precise how easy $\trop(g)$ is to work with algorithmically.  
In the theorem below, the underlying computational model is the {\em BSS 
model over $\R$} \cite{bcss}, and the {\em input size} of a 
point in $\Rn$ (resp.\ an $n$-variate $t$-nomial $g$) is 
defined to be $n$ (resp.\ $(n+1)t$), i.e., we merely measure the input size 
as the number of real numbers fed into a BSS machine.  
\begin{thm}
\label{thm:cxity2}
\scalebox{.97}[1]{Suppose $n$ is fixed. Then there is a polynomial-time
algorithm that, for any}\linebreak
input $r\!\in\!\Rn$ and $n$-variate exponential $t$-sum $g$,
outputs the closure --- described as an\linebreak
\scalebox{.98}[1]{explicit intersection of $O(t^2)$
half-spaces --- of the unique cell $\sigma_r$ of $\Sigma(\trop(g))$ containing
$r$.}
\end{thm}

\noindent 
We prove Theorem \ref{thm:cxity2} in Section \ref{sub:cxity2}. 

By applying the standard formula for point-hyperplane distance,
and the well-known efficient algorithms for approximating square-roots
(see, e.g., \cite{borwein}), Theorem \ref{thm:cxity2} implies that we can also 
efficiently check membership in any $\eps$-neighborhood about $\trop(g)$.  
Our complexity bound above, combined with our final main result below, tells us 
that membership in a neighborhood of $\trop(g)$ is a tractable and potentially 
useful relaxation of the problem of deciding membership in $\re(Z(g))$.
\begin{thm} \label{thm:intro2} 
Let $t\!\geq\!3$ and let $g$ be any $n$-variate exponential $t$-sum with 
spectrum $S\!:=\!\{a_1,\ldots,a_t\}\!\subset\!\Rn$, 
minimal frequency spacing 
$\delta(g)\!:=\!\min_{p\neq q} |a_p-a_q|$, and 
$d$ the dimension of the affine span of $S$. 
Then $d\!\leq\!\min\{n,t-1\}$ and:\\ 
\mbox{}\hspace{.95cm}(1) If $t\!=\!d+1$ then 
 $\trop(g)\!\subseteq\!\re(Z(g))$.\\ 
\mbox{}\hspace{.95cm}(2) If $t\!\geq\!d+1$ then (a) 
$\displaystyle{\sup
\limits_{\text{\scalebox{.7}[1]{$r\in\re(Z(g))$}}} \inf
\limits_{\substack{\mbox{} \\ 
\text{\scalebox{.7}[1]{$u\in\trop(g)$}}}} |r-u|\leq \frac{\log(t-1)}
{\delta(g)}}$ \ and\\ 
\mbox{}\hspace{4.55cm}(b) $\displaystyle{
\sup \limits_{\substack{\mbox{} \\ \text{\scalebox{.7}[1]{$u\in\trop(g)$}}}} 
\inf \limits_{\text{\scalebox{.7}[1]{$r\in\re(Z(g))$}}} 
|r-u|\leq \sqrt{ed}t^2\left.\left((\log 9)t
-\log\frac{81}{2}\right)\right/\delta(g)}$.\\  
(3) The bound from Assertion (2a) is optimal in the 
following sense: If $\delta\!>\!0$ and $\varphi(z)$ is\linebreak 
\mbox{}\hspace{.7cm}defined as $1+e^{\delta z_1}+\cdots+e^{\delta 
z_{t-1}}$ \ and \ $r\!:=\!-\log(t-1)
(1,\ldots,1)/\delta\!\in\!\R^{t-1}$, then 
$\re(Z(\varphi))\!\ni\!r$\linebreak 
\mbox{}\hspace{.7cm}and $\displaystyle{\inf
\limits_{\substack{\mbox{} \\ 
\text{\scalebox{.7}[1]{$u\in\trop(g)$}}}} |r-u|\!=\!(\log(t-1))/\delta}$.  
\end{thm}

\noindent 
\begin{minipage}[b]{0.5\linewidth}
\vspace{0pt}
\begin{ex}
When $g$ is the $2$-variate exponential $7$-sum
$\sum^6_{j=0}\binom{7}{j}e^{\cos(2\pi j/7)z_1
+\sin(2\pi j/7)z_2}$, Theorem \ref{thm:intro2} tells us that every point of 
$\re(Z(g))$ lies within distance\\
\mbox{$\log(6)/\sqrt{(1-\cos(2\pi/7))^2+\sin(2\pi/7)^2}
\!<\!2.065$}
of some point of $\trop(g)$. To the right, we can see $\trop(g)$ as the
black piecewise linear curve drawn on the right, along with the
\scalebox{.87}[1]{stated neighborhood of $\trop(g)$ containing $\re(Z(g))$.} 
\end{ex}
\end{minipage}\hspace{0cm} 
\begin{minipage}[b]{0.4\linewidth}
\vspace{0pt}
\mbox{\epsfig{file=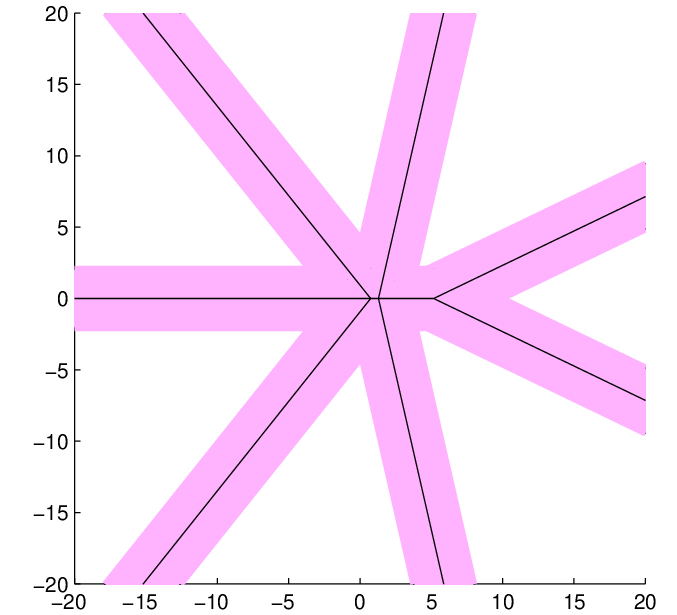,height=1.6in,clip=}\hspace{-.5cm}
\epsfig{file=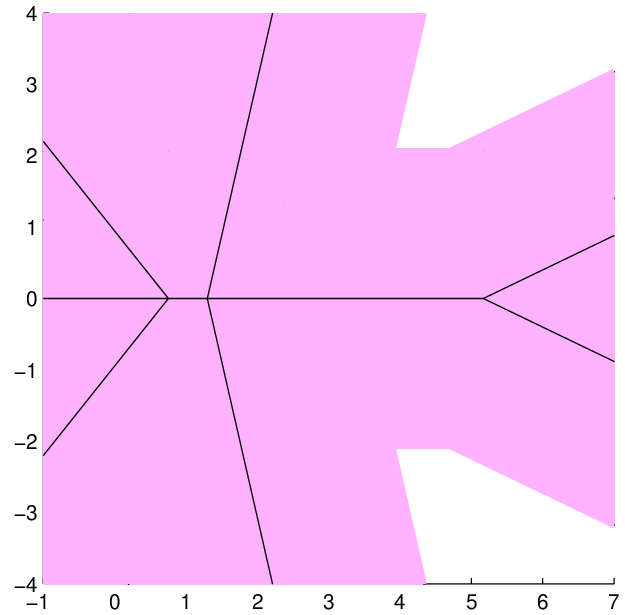,height=1.6in,clip=}}
\end{minipage}

\medskip 
\noindent 
We prove Theorem \ref{thm:intro2} in Section \ref{sec:proofs}. 
Prior to our work, there have been many fundamental results on the geometric 
and topological structure of the zero loci of exponential sums, e.g., 
\cite{moreno,kazarnovski,khovanski,favorov,silipo,soprunova,silipo2,
alessandrini,mora}. However, to the best of our knowledge, our results are 
the first to give an efficient approximation to all of $\re(Z(g))$ with 
explicit distance bounds.\footnote{A preliminary version of 
Theorem \ref{thm:intro2} appeared in our December 2014 Math ArXiV preprint 
{\tt 1412.4423 } and was presented by the first author at MEGA 2015 
(June 16, University of Trento).} 
Recently, Forsg\aa{}rd has found a bound 
complementary to Assertion (2a) of Theorem \ref{thm:intro2} that is tighter 
when the number of terms is exponential in the dimension. We rephrase his 
bound \cite[Thms.\ 1.2 \& 1.3]{forsgard} into our notation below:  
\begin{forthm} 
Following the notation of Theorem \ref{thm:intro2},\\ 
\mbox{}\hfill  
$\displaystyle{\sup
\limits_{\text{\scalebox{.7}[1]{$r\in\re(Z(g))$}}} \inf
\limits_{\substack{\mbox{} \\ 
\text{\scalebox{.7}[1]{$u\in\trop(g)$}}}} |r-u|\leq \frac{2n\sqrt{n}
\log(2+\sqrt{3})}{\delta(g)}}$.\hfill\mbox{}\\ 
In particular, if the spectrum of $g$ lies in $\Zn$, then the upper bound can be 
further improved to $n\log(2+\sqrt{3})$. \qed  
\end{forthm} 

\noindent 
For instance, for arbitrary real spectra, 
Forsg\aa{}rd's bound improves Assertion (2a) of our 
Theorem \ref{thm:intro2} when $t\!>\!1+e^{2.634n\sqrt{n}}$.   

One can also view the polyhedral structure in Theorem \ref{thm:intro2} as a 
limit shape of a parametric family of real parts of complex zero sets. 
Recall that, given any subsets $U,V\!\subseteq\!\Rn$, their
{\em Hausdorff distance} is 
$\Delta(U,V)\!:=\! 
\max\left\{\sup\limits_{u\in U}{} 
\inf\limits_{\substack{\mbox{}\\ v\in V}}|u-v|,
\sup\limits_{v \in V}{} \inf\limits_{\substack{\mbox{}\\ u\in U}}
|u-v|\right\}$. 
\begin{cor} \label{cor:limit}
For any exponential sum $g(z):=\sum_{j=1}^t e^{a_j\cdot z + \beta_j}$ we 
define a parametric\linebreak family of exponential sums via 
$g_s(z):=\sum_{j=1}^t e^{a_j\cdot z + s\cdot\beta_j}$ 
for any $s\!>\!0$. We then have\linebreak 
$\Delta\!\left( \frac{1}{s} \re(Z(g_s)) , \trop(g) \right)  \longrightarrow 0$ 
as $s \longrightarrow \infty$. 
\end{cor}

\noindent 
{\bf Proof:} 
First observe that $\trop(g_s)\!=\!s\trop(g)$ by definition. Applying Theorem 
\ref{thm:intro2} we then obtain $\Delta(\re(Z(g_s)),\trop(g_s))\leq 
\sqrt{en}t^2(2.2t-3.7)/\delta(g)$. So then,  
$\Delta(\re(Z(g_s)),\trop(g_s)) = s \Delta\!\left(\frac{1}{s}\re(Z(g_s)),
\trop(g)\right)$ and thus $\Delta\!\left(\frac{1}{s}\re(Z(g_s)),\trop(g)
\right)\leq \sqrt{en}t^2 (2.2t-3.7)/(\delta(g) s)$. \qed 

\medskip 
Corollary \ref{cor:limit} can be thought of as an exponential sum 
analogue of {\em Maslov dequantization}. The latter is a process by which 
one can obtain a (non-Archimedean) tropical variety as a limit of 
(complex) polynomial amoebae (see, e.g., \cite{virologpaper}).  

Let us now see a key ingredient, possibly of independent 
interest, behind the proof of our main multivariate metric bound. 

\subsection{Careful Projection to Reduce to the Univariate Case} 
Much of the recent literature on random projections aims toward creating
random matrices whose corresponding linear maps are ``nearly''
isometries. The approach is to create a random projection
matrix on a geometric object of interest, and the rank of the matrix is 
ultimately controlled by the statistical dimension of the geometric object
\cite{vershynin}. For our proof of Theorem \ref{thm:intro2}, we'll need a
projection of rank $1$ that distorts distances only slightly. Since most of 
the random matrix literature focusses on asymptotic behavior in high 
dimensions, we'll use a folkloric result stated as Lemma \ref{lemma:proj} 
below.

Let $G_{n,k}$ be the Grassmanian of $k$-dimensional subspaces of $\Rn$,
equipped with its unique rotation-invariant Haar probability measure
$\mu_{n,k}$. 
\begin{lemma} \label{lemma:proj} (See, e.g., \cite[Fact 3.2(c)]{milman2} 
and \cite[Lemma 6]{milman1}.)  
 Let $k\!\in\!\{1,\ldots,n-1\}$, $x\!\in\!\mathbb R^{n}$, and
$\varepsilon\!\leq\!\frac{1}{\sqrt{e}}$. Then
 \begin{equation*}
 \mu_{n,k}\! \left(\left\{ F\!\in\!G_{n,k} \; \left| \; 
|P_F(x)|\!\leq\!\varepsilon \sqrt{\frac{k}{n}}\cdot |x|\right.\right\}\right)
 \leq \left( \sqrt{e} \varepsilon\right)^{k}, 
 \end{equation*}
where $P_F$ is the surjective orthogonal projection mapping
$\Rn$ onto $F$. \qed
\end{lemma}

\noindent 
(See also \cite{vempala,vershynin} for more beautiful results on 
the theory and applications of random projections.) 
A simple consequence of Lemma \ref{lemma:proj} is the following existential
result.
\begin{prop}
\label{prop:proj}
Let $\gamma\!>\!0$ and  $x_{1}, \ldots , x_{N}\!\in\!\Rn$ be such
that $|x_{i}- x_{j}|\!\geq\!\gamma$ for all distinct $i,j$.
Then, following the notation of Lemma \ref{lemma:proj}, there is an  
$F\!\in\!G_{n,k}$ such that\linebreak
$\displaystyle{|P_F(x_{i}) - P_F(x_{j}) | \geq \sqrt{\frac{k}{en}}\cdot  
\frac{\gamma}{N^{2/k}}}$ for all distinct $i,j$.
\end{prop}

\noindent
{\bf Proof:}
Let $z_{\{i,j\}}\!:=\!|x_{i}- x_{j}|$. Then our assumption becomes
$z_{\{i,j\}}\!\geq\!\gamma$ for all distinct $i,j$, and there are
no more than  $N(N-1)/2$ such pairs $\{i,j\}$.
By Lemma \ref{lemma:proj} we have, for any fixed $\{i,j\}$, that
$|P_F(z_{\{i,j\}})|\!\geq\!\varepsilon \sqrt{\frac{k}{n}} z_{\{i,j\}}$ 
with probability at least $1-\left( \sqrt{e} \varepsilon\right)^{k}$. 
So the union bound for probabilities implies that,
for all distinct $i,j$, we have  
$|P_F(x_{i}) - P_F(x_{j})|\!\geq\!\varepsilon 
\gamma \sqrt{\frac{k}{n}}$ with probability at least
$1-\frac{N(N-1)}{2} (\sqrt{e} \varepsilon)^{k}$.
So our desired $F$ exists when 
$\varepsilon\!=\!\frac{1}{\sqrt{e} N^{2/k}}$ and we are done. \qed 

\medskip 
We now prove our main results.   

\section{Extending Classical Univariate Bounds to Exponential Sums: 
Proving Theorem \ref{thm:uni}} 
\label{sec:uni} 
The following simple quantitative bound on exponential sums will prove 
quite useful. In what follows, we let $[j]\!:=\!\{1,\ldots,j\}$. 
\begin{prop}
\label{prop:lop}
Suppose $t\!\geq\!3$ and $g(z_1)\!:=\!\sum^t_{j=1} e^{a_jz_1+\beta_j}$
satisfies $a_1\!<\cdots<\!a_t$ and $\beta_j\!\in\!\C$ for all $j$. Suppose
further that $u\!\in\!\trop(g)$, $\ell$ is the largest 
index such that\linebreak 
$\left|e^{a_\ell u+\beta_\ell}\right|\!=\!\max_{j\in[t]}  
\left|e^{a_j u+\beta_j}\right|$, 
and we set $\delta_\ell\!:=\!\!\!\!\!\!\!\!\min\limits_{\text{\scalebox{.7}[1]{
$p,q\in[\ell] \& p\neq q$}}} \!\!\!\!\!|a_p-a_q|$. 
Then for any $N\!\in\!\N$ and\linebreak  
$z_1\!\in\!\left[u+\frac{\log(N+1)}{\delta_\ell},\infty\right)\times 
\R$ we have $\left|\sum\limits^{\ell-1}_{j=1} e^{a_jz_1+\beta_j}\right|\!<\!
\frac{1}{N}\left|
e^{a_\ell z_1+\beta_\ell}\right|$.
\end{prop}  

\noindent
{\bf Proof:} First note that $2\!\leq\!\ell\!\leq\!t$ by construction. 
Let $b_j\!:=\!\re(\beta_j)$, $r\!:=\!\re(z_1)$, and observe  
\begin{eqnarray*}
\left|\sum^{\ell-1}_{j=1} e^{a_jz_1+\beta_j}\right| 
& \leq & \sum^{\ell-1}_{j=1} \left|e^{a_jz_1+\beta_j}\right| 
   = \sum^{\ell-1}_{j=1} e^{a_jr+b_j} 
   = \sum^{\ell-1}_{j=1} e^{a_j(r-u)+a_ju+b_j}.   
\end{eqnarray*}
Now, since $a_{j+1}-a_j\!\geq\!\delta_\ell$ for all $j\!\in\!\{1,
\ldots,\ell-1\}$, we obtain $a_j\!\leq\!a_\ell-(\ell-j)\delta_\ell$.
So for $r\!>\!u$ we have 
$\displaystyle{\left|\sum^{\ell-1}_{j=1} e^{a_jz_1+b_j}\right| 
\leq \sum^{\ell-1}_{j=1} e^{(a_\ell-(\ell-j)\delta_\ell)(r-u)+a_ju+b_j}\\ 
\leq \sum^{\ell-1}_{j=1} e^{(a_\ell-(\ell-j)\delta_\ell)(r-u)
                              +a_\ell u+b_\ell}}$,      
and thus 
\begin{eqnarray*}
\left|\sum^{\ell-1}_{j=1} e^{a_jz_1+b_j}\right| 
&\leq & e^{(a_\ell-(\ell-1)\delta_\ell)(r-u)+a_\ell u+b_\ell}
        \sum^{\ell-1}_{j=1} e^{(j-1)\delta_\ell (r-u)} \\ 
& = & e^{(a_\ell-(\ell-1)\delta_\ell)(r-u)+a_\ell u+b_\ell}
      \left(\frac{e^{(\ell-1)\delta_\ell(r-u)}-1}
      {e^{\delta_\ell(r-u)}-1}\right)\\
& < & e^{(a_\ell-(\ell-1)\delta_\ell)(r-u)+a_\ell u+b_\ell}
      \left(\frac{e^{(\ell-1)\delta_\ell(r-u)}}
      {e^{\delta_\ell(r-u)}-1}\right) =  \frac{e^{a_\ell r+b_\ell}}
      {e^{\delta_\ell(r-u)}-1}.  
\end{eqnarray*}
To prove our desired inequality, it thus clearly suffices to enforce
$e^{(r-u)\delta_\ell}-1\!\geq\!N$. The last inequality clearly
holds for all $r\!\geq\!u+\frac{\log (N+1)}{\delta_\ell}$, so we are done.
\qed

\medskip 
It is then easy to prove that the largest (resp.\ smallest) point of 
$\re(Z(g))$ can't be too much larger (resp.\ smaller) than the largest (resp.\ 
smallest) point of $\trop(g)$. Put another way, we can give an explicit 
vertical strip containing all the complex roots of $g$. 
\begin{cor} 
\label{cor:extremecauchy} Suppose $g$ is a univariate exponential $t$-sum  
with real spectrum and minimal frequency spacing 
$\delta(g)\!:=\!\min_{p\neq q} |a_p-a_q|$, 
$\um\!:=\!\min\trop(g)$, and $\up\!:=\!\max\trop(g)$.  
Then $\re(Z(g))$ is contained in the open interval 
$\left(\um-\frac{\log 2}{\delta(g)},\up+\frac{\log 2}{\delta(g)}
\right)$. 
\end{cor}  

\noindent 
Our earlier Lemma \ref{lemma:opt} tell us that the $\log 2$ in Corollary 
\ref{cor:extremecauchy} can not be replaced by 
any smaller constant. While 
the polynomial analogue of Corollary \ref{cor:extremecauchy} goes back 
to work of Cauchy, Birkhoff, and Fujiwara pre-dating 1916 
(see \cite[pp.\ 243--249, particularly bound 8.1.11 on pg.\ 247]{rs} and 
\cite{fujiwara} for further background)  
we were unable to find an explicit bound for exponential sums like Corollary 
\ref{cor:extremecauchy} in the literature. 
So we supply a proof below. 

\medskip 
\noindent 
{\bf Proof of Corollary \ref{cor:extremecauchy}:}  
Replacing $z_1$ by its negative, it clearly suffices to 
prove\linebreak 
$\re(Z(g))\!\subset\!\left(-\infty,\up+\frac{\log 2}{\delta(g)}\right)$. 
Writing $g(z_1)\!=\!\sum^t_{j=1}e^{a_jz_1+b_j}$ with $a_1\!<\cdots<\!a_t$, let 
$\zeta$\linebreak 
denote any root of $g$, $r\!:=\!\re(\zeta)$, and $\beta_j\!:=\!\re(b_j)$ for 
all $j$. Since we must have\linebreak 
$\sum^{t-1}_{j=1}e^{a_j\zeta+b_j}\!=\!-e^{a_t\zeta+b_t}$, taking 
absolute values implies that $\left|\sum^{t-1}_{j=1}e^{a_j\zeta+b_j}\right|
\!=\!\left|e^{a_t\zeta+b_t}\right|$.\linebreak  
\scalebox{.9}[1]{However, this equality is contradicted 
by Proposition \ref{prop:lop} for $\re(z_1)\!\geq\!\up+
\frac{\log 2}{\delta(g)}$. So we are done. \qed }   

\medskip  
Proposition \ref{prop:dumb} will then be a simple consequence of 
Corollary \ref{cor:extremecauchy} and the following special case 
of a fundamental result of Moreno. 
\begin{thm}
\label{thm:moreno} 
(Special case of \cite[Main Theorem, pg.\ 73]{moreno}.)  
Suppose $1,\alpha_1,\alpha_2,\alpha_3\!\in\!\R$ are linearly 
independent over $\Q$, $g(z_1)\!:=\!e^{\alpha_1 z_1}
+e^{\alpha_2 z_1} + e^{\alpha_3 z_1}$, $\sigma\!\in\!\R$, and 
the inequalities $|e^{\alpha_i \sigma}|\!\leq\!\sum_{j\in\{1,2,3\}\setminus
\{i\}} |e^{\alpha_j \sigma}|$ hold for all $i\!\in\!\{1,2,3\}$. Then 
$\sigma$ is a limit point of $\re(Z(g))$. \qed 
\end{thm} 

\medskip 
\noindent
{\bf Proof of Proposition \ref{prop:dumb}:} 
Let $g(z_1)\!:=\!e^{\sqrt{2}z_1}+e^{\sqrt{3}z_1}+e^{\sqrt{5}z_1}$. 
Clearly then, $\sqrt{3}-\sqrt{2}\!<\!\sqrt{5}-\sqrt{3}$, 
$\trop(g)\!=\!\{0\}$, and thus Corollary \ref{cor:extremecauchy} 
immediately implies the containment $X\!\subseteq\! 
\left(-\frac{\log 2}{\sqrt{3}-\sqrt{2}},\frac{\log 2}{\sqrt{3}-\sqrt{2}}
\right)$. Furthermore, since $g$ is an analytic function, its zeroes are 
isolated, and thus must be countable in number \cite{ahlfors}. 

Now note that 
$e^{\sqrt{5}u}\!>\!e^{\sqrt{3}u}\!>\!e^{\sqrt{2}u}$ for 
$u\!>\!0$, and this ordering is reversed for $u\!<\!0$. 
Furthermore, the same orderings apply to the corresponding 
derivatives. An elementary calculation then reveals that the hypothesis 
for Theorem \ref{thm:moreno} is satisfied at any $\sigma$ 
in the open interval $(-1.06,1.06)$. So we are done. \qed 

\medskip 
Our next result isolates vertical strips where no roots of $g$ can 
lie.  
\begin{cor} 
\label{cor:gap}  
Suppose $g(z_1)\!:=\!\sum^t_{j=1} e^{a_jz_1+\beta_j}$
satisfies $a_1\!<\cdots<\!a_t$, $\beta_j\!\in\!\C$ for all $j$, 
$\delta(g)\!:=\!\min_{p\neq q} |a_p-a_q|$, 
and that $u_1$ and $u_2$ are {\em consecutive} points of $\trop(g)$ 
satisfying $u_2\!\geq\!u_1+\frac{\log 9}{\delta(g)}$.  
Then the vertical strip 
$\left[u_1+\frac{\log 3}{\delta(g)},u_2-\frac{\log 3}{\delta(g)}\right]
\times \R$ contains {\em no} roots of $g$.  
\end{cor} 

\noindent 
{\bf Proof:} First note that $t\!\geq\!3$ since $\#\trop(g)\!\geq\!2$. 
Let $\ell$ be the unique index such that\\ 
\mbox{}\hfill $\left|e^{a_\ell u_1+\beta_\ell}\right|\!=\!\max_{j\in[t]}  
\left|e^{a_j u_1+\beta_j}\right|$ and
$\left|e^{a_\ell u_2+\beta_\ell}\right|\!=\!\max_{j\in[t]}  
\left|e^{a_j u_2+\beta_j}\right|$.\hfill\mbox{}\\  
There is a unique such index 
because, by the definition of $\trop(g)$, the point $(a_\ell,\re(\beta_\ell))$ 
lies at the intersection of two lines: One line goes through a  
pair of distinct points of the form $(a_i,\re(\beta_i))$ with 
$\left|e^{a_i u_1+\beta_i}\right|\!=\!\max_{j\in[t]}  
\left|e^{a_j u_1+\beta_j}\right|$, while the other 
goes through a pair of distinct points of the form $(a_k,\re(\beta_k))$ with  
$\left|e^{a_k u_2+\beta_k}\right|\!=\!\max_{j\in[t]}  
\left|e^{a_j u_2+\beta_j}\right|$.  

By Proposition \ref{prop:lop}, we have 
$\left|\sum^{\ell-1}_{j=1} e^{a_jz_1+\beta_j}\right|\!<\!\frac{1}{2}
 \left|e^{a_\ell z_1+\beta_\ell}\right|$ for all $z_1\!\in\!\left[u_1
+\frac{\log 3} {\delta(g)},\infty\right)$ and,\linebreak 
employing the change of variables $z_1\mapsto 
-z_1$, we obtain $\left|\sum^t_{j=\ell+1} e^{a_jz_1+\beta_j}
\right|\!<\!\frac{1}{2} \left|e^{a_\ell z_1+\beta_\ell}\right|$ for 
all\linebreak 
$z_1\!\in\!\left(-\infty,u_2-\frac{\log 3}{\delta(g)}
\right]$. 
So $\left|\sum_{j\neq \ell} e^{a_jz_1+\beta_j}\right|\!<\!
\left|e^{a_\ell z_1+\beta_\ell} \right|$ in the stated vertical strip, 
and this inequality clearly contradicts 
the existence of a root of $g$ in $\left[u_1+\frac{\log 3}{\delta(g)},u_2
-\frac{\log 3}{\delta(g)}\right]\times \R$. \qed 

\medskip 
An immediate consequence of Corollary \ref{cor:gap} is that the roots of $g$ 
always lie in the union of open vertical strips  
$\bigcup\limits_{u\in \trop(g)} \!\!\!\!\!\!\!\left(u-\frac{\log 3}
{\delta(g)}, u+\frac{\log 3}{\delta(g)}\right)\times \R$. It will in fact be 
the case that each connected component of this union contains roots of 
$g$ as well. To prove this, we will need some refined integral 
estimates. 

\subsection{Winding Numbers and Rectangles Around Tropical Points} 
It will be useful to first observe a basic fact on winding numbers along
line segments.
\begin{prop}
\label{prop:rouche}
Suppose $I\!\subset\!\C$ is any (compact) line segment and $g$ and $h$
are functions\linebreak
\scalebox{.885}[1]{analytic on a neighborhood of $I$ with
$|h(z)|<|g(z)|$ for all $z\!\in\!I$.
Then $\left|\im\left(\int_I\frac{g'+h'}{g+h}dz-\int_I\frac{g'}{g}dz\right)
\right|<\pi$.}
\end{prop}

\noindent
{\bf Proof:} The quantity $V_1\!:=\!\im\left(\int_I\frac{g'}{g}dz\right)$
(resp.\ $V_2\!:=\!\im\left(\int_I\frac{g'+h'}{g+h}dz\right)$)
is nothing more than the variation of the argument of $g$ (resp.\ $g+h$)
along the segment
$I$. Since $I$ is compact, $|g|$ and $|g+h|$ are bounded away from $0$ on $I$
by construction. So we can lift the paths $g(I)$ and $(g+h)(I)$ (in $\Cs$) to
the universal covering space induced by the extended logarithm function.
Clearly then, $V_1$ (resp.\ $V_2$) is simply a difference of values of
$\im(\Log(g))$ (resp.\linebreak $\im(\Log(g+h))$), evaluated at the endpoints
$I$, where different branches of $\Log$ may be used at each endpoint. In
particular, at any $z\!\in\!I$, our assumptions on $|g|$ and $|h|$
clearly imply that $g(z)+h(z)$ and $g(z)$ both lie in the open half-plane
normal (as a vector in $\R^2$) to $g(z)$. So $|\im(\Log(g(z)+h(z)))
-\im(\Log(g(z)))|\!<\!\frac{\pi}{2}$ at each of the two endpoints of $I$, 
and thus $|V_1-V_2|\!<\!\frac{\pi}{2}+\frac{\pi}{2}\!=\!\pi$. \qed

\medskip 
We will also need the following technical fact on the total variation of 
the imaginary part of $g'/g$ along horizontal line segments. 
\begin{thm} 
\label{thm:voorhoeve} \cite[Thm.\ 2]{voorhoevelms}  
Let $g(z_1)\!:=\!\sum^t_{j=1} e^{a_jz_1+\beta_j}$ with 
$a_1\!<\cdots<\!a_t$ and $\beta_j\!\in\!\C$ for all $j$. Also let 
$u,v\!\in\!\R$ with $g(u)g(v)\!\neq\!0$ and 
define $N$ to be the number of roots of $g$ on the closed interval 
$[u,v]$. Then  
$\displaystyle{\int^v_u \left|\im\!\left(\frac{g'(z)}{g(z)}\right)\right|dz
+N\pi \!\leq\!(t-1)\pi}$. \qed  
\end{thm} 

\noindent 
Note that since the $\beta_j$ are allowed to be complex, the bound above 
continues to hold if we integrate over any horizontal line segment in $\C$. 
Voorhoeve proved earlier in \cite[Lemma 1]{voorhoeveosc} that,  
for any $f$ meromorphic on an interval $[u,v]\!\subset\!\R$, the 
function $\im\!\left(\frac{f'}{f}\right)$ is analytic on $[u,v]$, 
save for a finite set of removable singularities. So the integral 
above is well-defined even if $g$ vanishes in the open interval $(u,v)$. 
\cite[Thm.\ 2]{voorhoevelms} in fact gives a sharper upper bound depending 
on the imaginary parts of the differences of the $\beta_j$, but we will only 
need the weaker bound stated above. See also \cite{voorhoeve} for an 
elegant and fascinating development of root counts for univariate 
exponential polynomials in various regions.    

\medskip 
We now state our final key root count behind Theorem \ref{thm:uni}. 
\begin{lemma} 
\label{lemma:vert} 
Let $g(z_1)\!:=\!\sum^t_{j=1} e^{a_jz_1+\beta_j}$ with $t\!\geq\!3$,  
$a_1\!<\cdots<\!a_t$, $\beta_j\!\in\!\C$ for all $j$, and let 
$\delta(g)\!:=\!\min_{p\neq q} |a_p-a_q|$, 
$\um\!:=\!\min\trop(g)$, and $\up\!:=\!\max\trop(g)$. 
Let $U_g$ be the union of open intervals $\left(\um-\frac{\log 2}{\delta(g)},\up+\frac{\log 2}{\delta(g)}
\right) \cap \bigcup\limits_{u\in \trop(g)} \!\!\!\!\!\!\!\left(u-\frac{\log 3}
{\delta(g)}, u+\frac{\log 3}{\delta(g)}\right)$. 
Let $\Gamma$ be any connected component of $U_g$ and let $p$ (resp.\ 
$q$) be the minimal (resp.\ maximal) index such that 
$\left|e^{a_p\cdot u+\beta_p}\right|\!=\!\max_j\left|e^{a_j\cdot u+\beta_j} 
\right|$ (resp.\ 
$\left|e^{a_q\cdot u+\beta_q}\right|\!=\!\max_j\left|e^{a_j\cdot u+\beta_j} 
\right|$) for some $u\!\in\!\Gamma$. 
Then $q\!>\!p$ and $g$ has at least one root in the rectangle 
$\Gamma\times \left[0,\frac{2(t+1)\pi}{\delta(g)}\right]$.  
\end{lemma} 

\medskip 
\noindent 
{\bf Proof of Lemma \ref{lemma:vert}:} That $q\!>\!p$ follows easily  
from the definition of $\trop(g)$: $\trop(g)\cap \Gamma$ is non-empty by 
construction, and if $u\!\in\!\Gamma\setminus\trop(g)$ 
then $\max_j\left|e^{a_j\cdot u+\beta_j}\right|$ is attained exactly once. 
Furthermore,  at least two terms of $g$ must be maximized in norm at any 
$u\!\in\!\trop(g)\cap \Gamma$, and $p$ 
(resp.\ $q$) must be no larger (resp.\ no smaller) than the index of any 
such term. 

Now let $\gm\!:=\!\inf \Gamma$ and $\gp\!:=\!\sup \Gamma$.  
Since $g$ is analytic, the Argument Principle (see, e.g., 
\cite{ahlfors}) 
tells us that the number of roots in our 
rectangle in question is exactly\\  
\mbox{}\hfill $\displaystyle{A:=\frac{1}{2\pi\sqrt{-1}}
\int_{I_-\cup I_+\cup J_-\cup J_+} \frac{g'(z)}{g(z)}dz}$\hfill\mbox{}\\ 
where $I_-$ (resp.\ $I_+$, $J_-$, $J_+$) is the 
oriented line segment from \\ 
\mbox{}\hfill $\left(\gm,\frac{2(t+1)\pi}{\delta(g)}\right)$ (resp.\ $(\gp,0)$, $(\gm,0)$,  
$\left(\gp,\frac{2(t+1)\pi}{\delta(g)}\right)$)\hfill\mbox{}\\ 
to\\ 
\mbox{}\hfill $(\gm,0)$ (resp.\ $\left(\gp,\frac{2(t+1)\pi}{\delta(g)}\right)$, 
$(\gp,0)$,  
$\left(\gm,\frac{2(t+1)\pi}{\delta(g)}\right)$),\hfill\mbox{}\\     
assuming no root of $g$ lies on $I_-\cup I_+\cup J_-\cup J_+$. 
By Corollaries \ref{cor:extremecauchy} and \ref{cor:gap}, there can 
be no roots of $g$ on $I_-\cup I_+$. So let assume temporarily that 
there are no roots of $g$ on $J_-\cup J_+$. 
 
By construction, any point of $\trop(g)\cap\Gamma$ is at least 
distance $\frac{\log 9}{\delta(g)}$ from any point of $\trop(g)\setminus 
\Gamma$. So Proposition \ref{prop:lop} tells us that when $p\!>\!1$ we have:\\  
\scalebox{.95}[1]{$\frac{1}{2}\left|e^{a_p \left(\gm+v\sqrt{-1} 
\right) +\beta_p}\right|\!>\!
\left|\sum\limits^{p-1}_{j=1}e^{a_j\left(\gm
+v\sqrt{-1} \right)+\beta_j} \right|$ 
\ \ and \ \  
$\frac{1}{2}\left|e^{a_p \left(\gm 
+v\sqrt{-1} \right) +\beta_p}\right|\!>\!
\left|\sum\limits^t_{j=p+1}e^{a_j\left(\gm +v\sqrt{-1} 
\right)+\beta_j} \right|$}\\ 
for all $v\!\in\!\R$. So then, 
$\left|e^{a_p \left(\gm+v\sqrt{-1} 
\right) +\beta_p}\right|\!>\!
\left|\sum\limits_{j\neq p}e^{a_j\left(\gm
+v\sqrt{-1} \right)+\beta_j} \right|$.   
(When $p\!=\!1$ Proposition \ref{prop:lop} 
yields the same conclusion in just one step.) 
So we can apply Proposition \ref{prop:rouche} and deduce that  
$\left|\im\left(\int_{I_-}\frac{g'(z)}{g(z)}dz-
\int_{I_-}\frac{(e^{a_p z+\beta_p})'}{e^{a_p z+\beta_p}}dz\right)
\right|\!<\!\pi$.  
So then, since $\int_{I_-}\frac{(e^{a_p z+\beta_p})'}
{e^{a_p z+\beta_p}}dz=\int_{I_-}a_p dz=
\frac{-2\pi\sqrt{-1}(t+1)a_p}{\delta(g)}$, we clearly obtain\\  
(1) \hfill 
$\displaystyle{\left|\im\left(\int_{I_-}
\frac{g'(z)}{g(z)}dz \right)-\frac{-2\pi\sqrt{-1}(t+1)a_p}{\delta(g)}
\right|<\pi.}$\hfill\mbox{}\\  
An almost identical argument (applying Propositions \ref{prop:lop} 
and \ref{prop:rouche} again, but with the term $\left|e^{a_q(\gp+v\sqrt{-1})
+\beta_q} \right|$ dominating instead) then yields\\ 
(2) \hfill  
$\displaystyle{\left|\im\left(\int_{I_+}
\frac{g'(z)}{g(z)}dz \right)-\frac{2\pi\sqrt{-1}(t+1)a_q}{\delta(g)}
\right|<\pi.}$\hfill\mbox{} 

So now we need only prove sufficiently sharp estimates on 
$\im\left(\int_{J_\pm}\frac{g'(z)}{g(z)}dz\right)$. Toward this end, observe that 
Theorem \ref{thm:voorhoeve} implies directly that 
$\displaystyle{\int_{J_\pm}\left|\im\!\left(\frac{g'(z)}{g(z)}\right)
\right|dz\!\leq\!(t-1)\pi}$. 
So combining with our estimates (1) and (2), 
and the additivity of integration, we obtain  
$\displaystyle{\left|A-\frac{(a_q-a_p)(t+1)}{\delta(g)}
\right|\!<\!t}$, in the special case where no roots of $g$ lie on $J_-\cup J_+$. 

To address the case where a root of $g$ lies on $J_-\cup J_+$, 
note that the analyticity of $g$ implies that the roots of $g$ are 
a discrete subset of $\C$. So we can find arbitrarily small $\eta\!>\!0$ with 
the boundary of the slightly stretched rectangle 
$\Gamma\times \left[-\eta,\frac{2(t+1)\pi}{\delta(g)}+\eta\right]$ 
not intersecting any roots of $g$, and define a similar normalized integral 
implementing the Argument Principle, 
which we'll call $A_\eta$, over the new contour. 
By the special case of our lemma already proved, we have 
$\displaystyle{\left|A_\eta-\frac{(a_q-a_p)\left(t+1+\frac{\eta\delta(g)}{\pi}\right)}{\delta(g)}
\right|\!<\!t}$. Let $n_\Gamma$ be the number of roots of 
$g$ in the rectangle $\Gamma\times \left[0,\frac{2\pi t}{\delta(g)}\right]$.  
Since $A_\eta\!=\!n_\Gamma$ for $\eta$ sufficiently small,  
we obtain $\displaystyle{\left|n_\Gamma-\frac{(a_q-a_p)(t+1)}{\delta(g)}
\right|\!\leq\!t}$. So\linebreak 
\scalebox{.92}[1]{$n_\Gamma\!\geq\!\frac{(a_q-a_p)(t+1)}{\delta(g)}-t\!\geq\!
\frac{(t+1)\delta(g)}{\delta(g)}-t\!=\!1$, and $g$ thus indeed has at 
least one root in $\Gamma\times \left[0,\frac{2(t+1)\pi}{\delta(g)}\right]$. 
\qed} 

\subsection{The Proof of Theorem \ref{thm:uni}:} 
\label{sub:uni} 
First note that the graph of $\cN_g$ is the lower hull of an 
intersection $P$ of exactly $t$ half-planes with edges of distinct slopes.  
So the polyhedron $P$ has at most $t$ edges, at most $t-1$ vertices, 
at least one vertex (since $t\!\geq\!2$), and 
thus the graph of $\cN_g$ has at most $t-1$ corners since corners correspond 
to vertices of $P$. We thus obtain $s\!\in\!\{1,\ldots,t-1\}$. 

\scalebox{.98}[1]{Assertion (1) on the containment $\re(Z(g))\!\subset\!U_g$ 
is immediate from Corollaries \ref{cor:extremecauchy} and \ref{cor:gap}.}  

Assertion (2) on each connected component of $U_g$ containing at least one 
point from $\re(Z(g))$ is immediate from Lemma \ref{lemma:vert}. 

To prove Assertion (3), we must show that near every point 
$u\!\in\!\trop(g)$ there is a point 
$r\!\in\!\re(Z(g))$ within distance $\frac{(\log 9)s-\log\frac{9}{2}}
{\delta(g)}$. (The remaining inequalities follow from the fact that 
$s\!\leq\!t-1$ and an elementary calculation.) 
So let $\Gamma$ be the unique connected component of $U_g$ containing 
$u$ and let $m\!:=\!\#(\trop(g)\cap \Gamma)$. We will prove that 
there is an $r\!\in\!\re(Z(g))$ within the distance 
$\frac{(\log 9)m-\log\frac{9}{2}}{\delta(g)}$ of $v$, yielding an  
inequality at least as tight as needed. Toward this end, observe first that 
consecutive points of $\trop(g)\cap \Gamma$ must be within distance 
strictly less than $\frac{\log 9}{\delta(g)}$. The maximal possible distance 
between $v$ and $r$ occurs when these two points lie at opposite extremes of 
the open interval $\Gamma$. Since $v$ must be no closer than 
$\frac{\log 2}{\delta(g)}$ to an endpoint of $\Gamma$, the maximal possible 
distance must be 
$\mathrm{Length}(\Gamma)-\frac{\log 2}{\delta(g)}\!<\!\frac{\log 2}{\delta(g)}
+\frac{\log 3}{\delta(g)}+\frac{\log 3}{\delta(g)}+(m-2)\frac{\log 9}{\delta(g)}
-\frac{\log 2}{\delta(g)}
\!=\!\frac{(\log 9)m-\log\frac{9}{2}}{\delta(g)}$, assuming (without loss of 
generality) that $v$ is as far left as possible and $r$ is as far right as 
possible. \qed 

\section{Algorithmic Complexity: The Proofs of Theorems \ref{thm:mem} 
and \ref{thm:cxity2}}
\label{sec:cxity} 
The {\em BSS model over $\R$} \cite{bcss} naturally augments the classical 
{\em Turing machine} \cite{papa,arora,sipser2} by    
allowing field operations and comparisons over $\R$ in unit time. 
We are in fact forced to move beyond the Turing model since our 
exponential sums involve arbitrary real numbers,  
and the Turing model only allows finite bit strings as inputs.  
We refer the reader to \cite{bcss} for further background. 

We recall here some basic facts about the set of inputs on which a BSS machine 
over $\R$ terminates.  
\begin{thm} 
\label{thm:bcss} \cite[Thm.\ 1, Pg.\ 52]{bcss} 
The halting set of a BSS machine over $\R$ is a countable 
union of semi-algebraic sets. \qed 
\end{thm} 
\noindent
The converse of Theorem \ref{thm:bcss} fails in general:  
For instance, if $S$ is any countably infinite subset of
a transcendence basis for $\R$ over $\Q$, then $S$ can not
be the halting set of any BSS machine over $\R$. (One can
even write such subsets in terms of infinite series,
via an explicit basis found by von Neumann \cite{vonneumann} around
1928.) This follows immediately from the following consequence of the
development in \cite[Sec.\ 2.3]{bcss}:
\begin{sop}
\label{prop:count}
Any countable subset of $\R$ that is the halting set for a BSS machine over
$\R$ must be a subset of the algebraic closure of a real extension of $\Q$ of
finite transcendence degree. \qed
\end{sop}

Let us also recall the following basic facts about {\em semi-algebraic} sets,
i.e., the solution sets of finite collections of polynomial inequalities
and polynomial equalities in $\Rn$: First, semi-algebraic sets are
closed under all Boolean operations (intersection, union, and complement).
Also, semi-algebraic sets admit a natural notion of dimension, via
the largest $d$ permitting a semi-algebraic embedding of a real $d$-ball
(see, e.g., \cite[Ch.\ 5, Sec.\ 5.3, pp.\ 170--172]{bpr}). Some additional
qualitative facts we'll also need can be summarized
as follows:
\begin{tame}
Suppose $S\!\subset\!\R^2$ is semi-algebraic, and $\bar{S}$ and
$S^{\circ}$ respectively denote the closure and interior of $S$. Then: \\
1.\  $\bar{S}$, $S^\circ$, and $\partial S\!:=\!\bar{S}\setminus S^\circ$ are semi-algebraic.\\
2.\ $S$ has only finitely many connected components, each of
which is semi-algebraic.\\
3.\ If, in addition, $S$ is a connected curve, then $S$ has only finitely
many singularities. \\
4.\ Let $\rho : \R^2\longrightarrow \R$ denote the
projection defined by $\rho(x,y)\!=\!x$. Then, continuing
Assertion\linebreak
\mbox{}\hspace{.5cm}(3), there is an $n_S\!\in\!\N$ such that all fibers of
$\rho$ have cardinality at most $n_S$.
\end{tame}

\noindent 
{\bf Note:} {\em Neither Proposition BCS nor the preceding tameness theorem 
appeared in the published (Mathematische Annalen, Vol.\ 377, pp.\  
863--882 (2020)) version of this paper. In particular, we 
inserted the theorem above after Alexander Rashkovskii kindly pointed out in 
late July 2020 that our published proof of Theorem \ref{thm:mem} had an error. 
We inserted Proposition BCS after a discussion with Lenore Blum, Felipe Cucker, 
and Mike Shub, on the simplest failures of the converse of Theorem 3.1. 
We apply the tameness theorem to give a corrected proof of Theorem 
\ref{thm:mem} below, but we will have no further need for Proposition BCS. 
\dia} 

\medskip 
\noindent
{\bf Proof of the Semi-Algebraic Tameness Theorem:} The first two assertions of 
(1) are exactly the content of
\cite[Ch.\ 3, Prop.\ 3.1, pg.\ 84]{bpr}. The final assertion of (1) is then
immediate since semi-algebraic sets are closed under Boolean operations
by definition.

Assertion (2) is immediate from the notion of
{\em cylindrical decomposition}. The latter is a refined decomposition of
a semi-algebraic set into finitely many (semi-algebraic) connected components,
and the existence of such a decomposition is a classical fact: See, e.g.,
\cite[Ch.\ 5, Thm.\ 5.6, pg.\ 163]{bpr}. In particular, the tameness
of fibers from Assertion (4) is also an immediate consequence of
cylindrical decomposition.

Assertion (3) is a direct consequence of the notion of {\em semi-algebraic
cell stratification of $\R^2$ adapted to $S$}. The latter is a partition $S$
into finitely many semi-algebraic smooth manifolds (here, each diffeomorphic
to an open interval or a point) called {\em strata}, such that the
closure of any stratum is a union of strata. That such stratifications
exist (and in much greater generality) is also a classical fact: See, e.g.,
\cite[Ch.\ 5, Thm.\ 5.38, pg.\ 177]{bpr}. \qed

\subsection{The Proof of Theorem \ref{thm:mem}} 
\label{sub:mem}  
Let $\Cs\!:=\!\C\setminus\{0\}$ and let\\
\mbox{}\hspace{.5cm}
$R\!:=\!\re(Z(1-e^{z_1}-e^{z_2}))$ \ and \
$S\!:=\!\{(\log|x|,\log|y|)\; | \; 1-x-y=0; \ 
x,y\!\in\!\Cs\}$. \\
Via the equality $\log\left|e^{\alpha+\sqrt{-1}\beta}\right|
\!=\!\alpha$ (valid for any $\alpha,\beta\!\in\!\R$)
we see that\\
\scalebox{.89}[1]{$(x,y)\!\in\!\re(Z(1-e^{z_1}-e^{z_2})) \Longleftrightarrow 
(x,y)\!=\!(\log|e^{z_1}|,\log|e^{z_2}|)$ for some $(z_1,z_2)\!\in\!\C^2$ with
$1-e^{z_1}-e^{z_2}\!=\!0$.}\linebreak
Since the exponential function defines a surjection from $\C$ onto $\Cs$ we
then clearly have $R\!=\!S$.

Now note that $J\!:=\!\{(|w_1|,|w_2|)\; | \; 1-w_1-w_2\!=\!0; \ w_1,
w_2\!\in\!\Cs\}$ is exactly the following\linebreak
\scalebox{.94}[1]{semi-infinite strip with corners
deleted: $I\!:=\!\{(x,y)\!\in\!\R^2\; | \; -1\!\leq\!y-x\!\leq\!1, \  
x+y\!\geq\!1, \text{ and } xy\!\neq\!0\}$.}\linebreak
This is because $w_1+w_2\!=\!1 \Longrightarrow |w_1+w_2|\!=\!1$,
$|w_1|\!=\!|1-w_2|\!=\!|w_2-1|$, and $|w_2|\!=\!|1-w_1|\!=\!|w_1-1|$. So by the
Triangle Inequality we obtain $|w_1|+|w_2|\!\geq\!1$, $|w_1|\geq\!||w_2|-1|$,
and $|w_2|\geq\!||w_1|-1|$, and thus (setting $x\!=\!|w_1|$ and $y\!=\!|w_2|$)
we obtain $J\!\subseteq\!I$. To see that
$I\!\subseteq\!J$, assume $(x,y)\!\in\!I$ and consider
$y_\theta\!:=\!1+xe^{\theta\sqrt{-1}}$ for $\theta\!\in\![0,\pi]$.
Clearly $x\!>\!0$. So then $|y_\theta|^2\!=\!(1+(\cos \theta)x)^2+
(\sin \theta)^2x^2\!=\!1+2(\cos \theta)x+x^2$ is a decreasing
differentiable function of $\theta$,
with $|y_0|\!=\!x+1$ and $|y_\pi|\!=\!|x-1|$. Since $|x-1|\!\leq\!y\!\leq\!x+1$
there must then be a $\theta\!\in\![0,\pi]$ with $y\!=\!|y_\theta|$.
Letting $w_1\!:=\!-xe^{\theta\sqrt{-1}}$ and $w_2\!:=\!y_\theta$, we then
obtain $w_1+w_2\!=\!1$, $|w_1|\!=\!|x|\!=\!x$, and $|w_2|\!=\!y$. So we have
obtained $I\!\subseteq\!J$ and thus $I\!=\!J$.

Clearly then, $R$ is simply the image of $I$
under the (differentiable) coordinate-wise logarithm map. In particular, we
see that the curve $Y$ defined by $y\!=\!\log(1+e^x)$, as $x$ ranges over all
of $\R$, is a connected component of the boundary $\partial R$.

By Theorem \ref{thm:bcss}, if membership in $R$ is decidable, then
$R$ must be a countable union $\bigcup_{i\in \N} S_i$ of semi-algebraic sets
$S_i$. Let $W\!:=\!Y\cap([0,1]\times \R)$, abusing notation slightly by
identifying $\C$ with $\R^2$. Then
$W$ is compact and infinite, and thus some $S_i$ must have $W\cap S_i$
infinite. Note in particular that $W\cap S^\circ_i\!=\!\emptyset$ (since
$S^\circ_i$ is in the interior of $R$) and thus (by the Semi-Algebraic 
Tameness Theorem) $S_i\setminus S^\circ_i$ must be a finite union of isolated
points and smooth connected semi-algebraic curves.
In particular, $S_i$ must contain a smooth connected semi-algebraic curve $C$
such that $W\cap C$ is infinite. Recalling that $\rho : 
\R^2\longrightarrow \R$ is the projection defined by $\rho(x,y)\!=\!x$,
we may assume further that $C$ is the graph of a smooth
algebraic function $f$ on a non-empty open sub-interval of $(0,1)$,
via the Implicit Function Theorem and Assertion (4) of the Semi-Algebraic 
Tameness Theorem.  (In particular, this might entail replacing $C$
with a non-empty, connected (and semi-algebraic), open subset of $C$.)

Now observe that $W\cap C$ (resp.\ $\rho(W\cap C)$) must have at least one
accumulation point since $W$ (resp.\ $[0,1]$) is compact, and thus
the graphs of the smooth functions $\log(1+e^x)$ and $f$ agree on
an infinite sequence of points with a limit point.
But this is impossible, since $\log(1+e^x)$ is a transcendental function.
In particular, since $\log(1+e^x)$ is analytic on the
domain $\R\times (-\pi,\pi)$, the function $f$ must have an analytic
continuation to an algebraic function with an essential singularity at
$\infty$ \cite[Pg.\ 127]{ahlfors}. Since algebraic functions can only have
zeroes or poles of finite fractional order at $\infty$, we obtain a
contradiction.
\qed

\medskip 
\noindent 
{\bf Note:} {\em The key obstruction to membership in $\re(Z(g))$ being  
decidable --- demonstrated above --- is that the {\em boundary} of 
$\re(Z(1-e^{z_1}-e^{z_2}))$ is not expressible as a countable union of 
semi-algebraic sets. Indeed, it is easy to express the interior of 
$\re(Z(1-e^{z_1}-e^{z_2}))$ as a countable union of disks: simply 
consider the union of all open disks contained in\linebreak  
\scalebox{.95}[1]{$\re(Z(1-e^{z_1}-e^{z_2}))$, centered at some rational 
point in the interior of $\re(Z(1-e^{z_1}-e^{z_2}))$. \dia}}

\subsection{Proving Theorem \ref{thm:cxity2}}  
\label{sub:cxity2} 
We will need some supporting results on linear programming before starting our 
proof of Theorem \ref{thm:cxity2}. The results we'll need are covered with 
great clarity in well-known monographs such as \cite{schrijver,grotschel,
gritzmann}.
\begin{dfn} 
Given any matrix $M\!\in\!\R^{N\times n}$ with $i\thth$ row $m_i$, and
$c\!:=\!(c_1,\ldots,c_N)\!\in\!\R^N$, the notation $Mx\!\leq\!c$ means
that $m_1\cdot x\!\leq\!c_1,\ldots,m_N\cdot x\!\leq\!c_N$ all hold. 
These inequalities are called {\em constraints}, and 
the set of all $x\!\in\!\R^N$ satisfying $Mx\!\leq\!c$ is called the 
{\em feasible region} of $Mx\!\leq\!c$. We also call  
a constraint {\em active} if and only if it holds with equality. 
Finally, we call a constraint {\em redundant} if and only if
the corresponding row of $M$ and corresponding entry of $c$ can be deleted
without affecting the feasible region of $Mx\!\leq\!c$. \dia 
\end{dfn} 
\begin{lemma}
\label{lemma:red} 
Suppose $n$ is fixed. Then, given any $c\!\in\!\R^N$ and 
$M\!\in\!\R^{N\times n}$, we can, in time polynomial in $N$,
find a submatrix $M'$ of $M$, and a subvector $c'$ of $c$,
such that the feasible regions of $Mx\!\leq\!c$ and $M'x\!\leq\!c'$ are 
equal, and $M'x\!\leq\!c'$ has no redundant constraints. 
Furthermore, in time polynomial in $N$, we can also enumerate all 
maximal sets of active constraints defining vertices of the 
feasible region of $Mx\!\leq\!c$. \qed
\end{lemma} 

\noindent 
Note that we are using the BSS model over $\R$ in the preceding 
lemma. In particular, we are counting just field operations and 
comparisons over $\R$, and these are the only operations needed 
above. There are many possible choices for the underlying algorithm: 
For instance, the classical Simplex Algorithm (using, say, Bland's 
Anticycling Rule) very easily yields Lemma 
\ref{lemma:red}. Note that the assumption that 
$n$ be fixed is critical: As of October 2018, it is still an open problem 
whether Linear Programming can be done in time polynomial in {\em both} 
$n$ and $N$ in the BSS model over $\R$ (a.k.a.\ {\em strongly polynomial-time}  
in older terminology).  

\medskip 
\noindent 
{\bf Proof of Theorem \ref{thm:cxity2}:}  
Let $r\!\in\!\Rn$ be our input query point. 
Let $b_j\!:=\!\re(\beta_j)$ for all $j$. 
Using $O(t\log t)$ comparisons and $O(n)$ arithmetic operations, 
we can first isolate all indices such that
$\max_j\{a_j\cdot r+b_j\}$ is attained, so let $j_0$ be any such 
index. (Note that these are the same indices we would obtain if we 
were to maximize $\left|e^{a_j\cdot z+\beta_j}\right|$.) We then obtain, say, 
$J$ equations of the form $a_j\cdot r+b_j\!=\!a_{j_0}\cdot r+b_{j_0}$ 
and $K$ inequalities of the form 
$a_j\cdot r+b_j\!<\!a_{j_0}\cdot r+b_{j_0}$. 

Thanks to Lemma \ref{lemma:red}, we can determine the exact cell of $\trop(f)$ 
containing $r$ if $J\!\geq\!2$. Otherwise, we obtain the unique cell of 
$\Rn\!\setminus\!\trop(f)$ with relative interior containing $r$. Note also 
that an $(n-1)$-dimensional face of either kind of cell must be contained 
in a hyperplane of the form $\{u\!\in\!\Rn\; | \; (a_{j_1}-a_{j_2})\cdot u 
+(b_{j_1}-b_{j_2})\!=\!0\}$ for some distinct indices $j_1$ and $j_2$. 
So there are at most $t(t-1)/2$ such $(n-1)$-dimensional faces, and thus 
$\sigma_r$ is the intersection of at most $t(t-1)/2$ half-spaces. So we are 
done. \qed

\section{The Proof of Our Main Multivariate Bound: Theorem \ref{thm:intro2}}
\label{sec:proofs} 
Let us first observe that $d\!\leq\!\min\{n,t-1\}$ follows immediately 
from the basic fact that any $d$-polytope in $\Rn$ has dimension at most 
$n$ and at least $d+1$ vertices. 

In what follows, for any real $n\times n$ matrix $M$ and 
$z\!\in\!\Rn$, we assume that $z$ is a column vector when we write $Mz$.
Also, for any subset $S\!\subseteq\!\Rn$, the notation $MS\!:=\!\{Mz\; | \; 
z\!\in\!S\}$ is understood. The following simple functoriality properties of
$\trop(g)$ and $\re(Z(g))$ will prove useful.
\begin{prop}
\label{prop:rescale}
Suppose $g_1$ and $g_2$ are $n$-variate exponential $t$-sums,
$\alpha\!\in\!\Cs$, $a\!\in\!\Rn$,
$\beta\!:=\!(\beta_1,\ldots,\beta_n)\!\in\!\Cn$,
and $g_2$ satisfies the identity $g_2(z)\!=\!\alpha e^{a\cdot z} 
g_1(z_1+\beta_1,\ldots,z_n+\beta_n)$. Then
$\re(Z(g_2))\!=\!\re(Z(g_1))-\re(\beta)$ and
$\trop(g_2)\!=\!\trop(g_1)-\re(\beta)$. Also, if $M\!\in\!\R^{n\times n}$
and we\linebreak
\scalebox{.89}[1]{instead have the identity $g_2(z)\!=\!g_1(Mz)$, then
$M\re(Z(g_2))\!=\!\re(Z(g_1))$ and $M\trop(g_2)\!=\!\trop(g_1)$. \qed}
\end{prop}

\subsection{The Proof of Assertion (1) of Theorem \ref{thm:intro2}} 
First note that, thanks to Proposition \ref{prop:rescale}, an invertible linear 
change of variables allows us to reduce to the special case 
$\{ a_1, \ldots, a_{n+1} \}=\!\{\bO,e_1,\ldots,e_n\}$,  
where $\bO$ and $\{e_1,\ldots,e_n\}$ are respectively the origin 
and standard basis vectors of $\Rn$. But this special case is 
well-known: One can either prove it directly, or avail to earlier 
work on the spines of amoebae, e.g., \cite[Thms.\ 1 \& 2]{passarerullgard}. 
(See also  \cite[Prop.\ 3.1.8]{forsbergthesis} or the remark following 
Theorem 8 on Page 33, and Theorem 12 on Page 36, of \cite{rullgardthesis} 
for precursors.)   
\qed 

\subsection{The Proof of Assertion (2a) of Theorem \ref{thm:intro2}}  
By Assertion (2b) (proved independently in the 
next section) $Z(g)$ is non-empty. So pick any $z\!\in\!Z(g)$, 
let $r\!:=\!\re(z)$, and assume without 
loss of generality that $\left|e^{a_1\cdot z+\beta_1}\right|\!\geq\!
\left|e^{a_2\cdot z+\beta_2}\right|
\!\geq \cdots \geq\!\left|e^{a_t\cdot z+\beta_t}\right|$. 
Since $g(z)\!=\!0$\linebreak implies  
$\left|e^{a_1\cdot z+\beta_1}\right|\!=\!\left|e^{a_2\cdot z+\beta_2}+
\cdots+e^{a_t\cdot z+\beta_t}\right|$, the Triangle Inequality 
immediately implies that\linebreak  
$\left|e^{a_1\cdot z+\beta_1}\right|\!\leq\!(t-1)\left|e^{a_2\cdot z+
\beta_2} \right|$. Letting $b_j\!:=\!\re(\beta_j)$ for all $j$ and then taking 
logarithms we obtain 
\setcounter{equation}{2} 
\begin{eqnarray} 
\label{eqn:3} 
& a_1\cdot r+b_1\geq \cdots \geq a_t\cdot r +b_t& \text{ and}  
\end{eqnarray} 
\begin{eqnarray} 
\label{eqn:4} 
& a_1\cdot r +b_1\leq \log(t-1)+a_2\cdot r +b_2& 
\end{eqnarray} 
For each $j\!\in\!\{2,\ldots,t\}$ let us then define $\eta_j$ to be
the shortest vector such that \\ 
\mbox{}\hfill $a_1\cdot (r+\eta_j)+b_1 
=a_j\cdot (r+\eta_j)+b_j$. \hfill\mbox{}\\ 
Note that $\eta_j\!=\!\lambda_j(a_j-a_1)$ for some nonnegative $\lambda_j$ 
since we are trying to affect the dot-product $\eta_j\cdot(a_1-a_j)$. 
In particular, 
$\lambda_j\!=\!\frac{(a_1-a_j)\cdot r+b_1-b_j}{|a_1-a_j|^2}$ so
that $|\eta_j|\!=\!\frac{|(a_1-a_j)\cdot r + b_1-b_j|}
{|a_1-a_j|}$. (Indeed, Inequality (\ref{eqn:3})   
implies that $(a_1-a_j)\cdot r+b_1-b_j\!\geq\!0$.)  

Inequality (\ref{eqn:4}) implies that 
$(a_1-a_2)\cdot r+b_1-b_2\!\leq\!\log(t-1)$. 
We thus obtain\linebreak $|\eta_2|\!\leq\!\frac{\log(t-1)}
{|a_1-a_2|}\!\leq\!\frac{\log(t-1)}{\delta(g)}$.
So let $j_0\!\in\!\{2,\ldots,t\}$ be any $j$ minimizing
$|\eta_j|$. We of course have $|\eta_{j_0}|\!\leq\!(\log(t-1))/\delta(g)$ and,  
by the definition of $\eta_{j_0}$, we have\\ 
\mbox{}\hfill  
$a_1\cdot (r+\eta_{j_0})+b_1\!=\!a_{j_0}\cdot (r+\eta_{j_0})
+b_{j_0}$.\hfill\mbox{}\\ 
Moreover, the fact that $\eta_{j_0}$ is the shortest among the $\eta_j$ 
implies that\\
\mbox{}\hfill $a_1\cdot (r+\eta_{j_0})+b_1\!\geq\!a_j\cdot (r+ 
\eta_{j_0})+b_j$\hfill\mbox{}\\  
for all $j$: Otherwise, for some $j'$, we would have $a_1\cdot 
(r+\eta_{j_0})+b_1\!<\!a_{j'}\cdot (r+\eta_{j_0})+b_{j'}$ and
$a_1\cdot r+b_1\!\geq\!a_{j'}\cdot r+b_{j'}$ (the latter 
following from Inequality (\ref{eqn:3})). Taking a convex linear 
combination of the last two inequalities, we would then obtain  
a $\mu\!\in\![0,1)$ such that\\ 
\mbox{}\hfill $a_1\cdot (r+\mu\eta_{j_0})+b_1\!=\!a_{j'}\cdot 
(r+\mu\eta_{j_0})+b_{j'}$.\hfill\mbox{}\\ 
Thus, by the definition of $\eta_{j'}$, we would obtain $|\eta_{j'}|
\!\leq\!\mu|\eta_{j_0}|\!<\!|\eta_{j_0}|$ --- a contradiction.

We thus have (i) $a_1\cdot (r+\eta_{j_0})+b_1\!=\!
a_{j_0}\cdot (r+\eta_{j_0})+b_{j_0}$, 
(ii) $a_1\cdot (r+\eta_{j_0})+b_1\!\geq\!
a_j\cdot (r+\eta_{j_0})+b_j$ for all $j$, and (iii) 
$|\eta_{j_0}|\!\leq\!(\log (t-1))/\delta(g)$. 
Together, these inequalities imply that 
$u\!:=\!r+\eta_{j_0}\!\in\!\trop(g)$ and $|r-u|\!\leq\!(\log(t-1))/
\delta(g)$. \qed 

\subsection{The Proof of Assertion (2b) of Theorem \ref{thm:intro2}}
Thanks to Proposition \ref{prop:rescale}, we can apply a suitable 
orthogonal linear change of variables to assume that  
$d\!=\!n$. By the $k\!=\!1$ case of Proposition 
\ref{prop:proj} we then deduce that there 
exists a unit vector $\theta\!\in\!\Rn$ such that 
\begin{gather}\label{proj}
\min_{i\neq j}|a_i\cdot \theta - a_j\cdot \theta | \geq 
\frac{\delta(g)}{\sqrt{en} t^2} 
\end{gather}
Now let $u\!\in\!\trop(g)$ and write 
$u\!=\!u_\theta \theta + u^\bot_\theta$ for some $u_\theta\!\in\!\R$  
and $u^\bot_\theta\!\in\!\Rn$ 
perpendicular to $\theta$.  Our goal is to find $z\!\in\!\Cn$ with 
$g(z)\!=\!0$ and $|\re(z)-u|\!\leq\!\frac{\sqrt{en}t^2
\left((\log 9)t-\log\frac{81}{2}\right)}{\delta(g)}$.

For $z_1\!\in\!\C$ we then define the univariate exponential $t$-sum  
$\g(z_1)\!:=\!\sum_{j=1}^{t} e^{a_j \cdot z_1 \theta + a_j \cdot 
u^\bot_\theta + \beta_j}$. $\g$ is in fact the restriction of $g$ to the 
complex line parametrized by $l(z_1)=z_1 \theta +  u^\bot_\theta$. In 
particular, $\g$ has the same number of terms as $g$ thanks to our choice of 
$\theta$,  and the definition of $\trop(\g)$ implies that 
$u_\theta\!\in\!\trop(\g)$. By Theorem \ref{thm:uni} there is an 
$\omega\!\in\!\mathbb{C}$ such that 
$0\!=\!\g(\omega)\!=\!g(l(\omega))$ and $|\re(\omega) - u_{\theta}| 
\leq \frac{\left((\log 9)t-\log\frac{81}{2}\right)}{\delta\!\left(\g\right)}$.  
Since $|\re(l(\omega))- u|\!=\!|(\re(\omega)-u_{\theta})\theta|$, 
and $\delta\!\left(\g \right)\!\geq\!\frac{\delta(g)} {\sqrt{en} t^2}$ 
by Inequality (\ref{proj}), we can conclude by taking $z\!:=\!l(\omega)$. \qed 

\subsection{The Proof of Assertion (3) of Theorem \ref{thm:intro2}}
Since $\varphi(r)\!=\!0$ it is clear that $\re(Z(g))\!\ni\!r$. 
It is easily checked that $\trop(g)$ is the codimension $1$ part of 
the outer normal fan of the standard $n$-simplex in $\Rn$. So 
$r$ is in fact at distance $(\log(t-1))/\delta$ from $\trop(g)$ 
because $r$ lies in the negative orthant and is at distance 
$(\log(t-1))/\delta$ from each of the coordinate hyperplanes of 
$\Rn$. \qed 

\section*{Acknowledgements} 
We thank  
Timo de Wolff for pointing out Silipo's 
work \cite{silipo}. We also thank Pascal Koiran, Gregorio Malajovich,  
Ji\v{r}\'{\i} Matou\v{s}ek, and Klaus Meer for useful discussions. 
Special thanks go to our referee for valuable commentary that significantly 
improved our paper.  

\medskip 
\noindent 
{\bf Note:} {\em We also thank Alexander Rashkovskii who, in late July 2020, 
pointed out an error in the published proof of Theorem 1.1. We have corrected 
the proof of Theorem 1.1 in this version of our paper. We also thank Lenore 
Blum, Felipe Cucker, and Mike Shub for discussions leading us to Proposition 
BCS. \dia} 

\medskip 
In closing, we would like to remember our friend and colleague Joel Zinn: 
He was admired at Texas A\&{}M (and far beyond) for his wisdom, warmth, humor, 
and kindness. He is sorely missed. 

\bibliographystyle{acm}

\end{document}